# Solutions of systems of the partial differential equations of Kampé de Fériet type functions

Anvar Hasanov[1] and Tuhtasin Ergashev[2]

[1] Institute of Mathematics, 29 Durmon yuli street, Tashkent 100125, Uzbekistan
E-mail: anvarhasanov@yahoo.com

[2] The Tashkent institute of engineers of irrigation and mechanization of agriculture. 39 Kari-Niyazi street, Tashkent, Uzbekistan. E-mail: ertuhtasin@mail.ru

**Abstract**

In investigation of boundary-value problems for certain partial differential equations arising in applied mathematics we often need to study the solution of system of partial differential equations satisfied by hypergeometric functions and find explicit linearly independent solutions for the system. In this present investigation, we give the solutions of systems of partial differential equations for two Kampé de Fériet type functions of third and fourth orders and of two variables.

## 1. Introduction

A great interest in the theory of hypergeometric functions (that is, hypergeometric functions of one, two and several variables) is motivated essentially by the fact that solutions of many applied problems involving thermal conductivity and dynamics, electromagnetic oscillation and aerodynamics, and quantum mechanics and potential theory are obtainable with the help of hypergeometric (higher and special or transcendent) functions [5,8,17,18]. Such kinds of functions are often referred to as special functions of mathematical physics. For the purpose of the present work, we recall the following definition of the most general hypergeometric function of two variables $F_{l:i;j}^{p:q;k}$, that is the Kampé de Fériet hypergeometric series of two variables (see [19]):

$$F_{l:i;j}^{p:q;k}\left[\begin{matrix}(a_p):(b_q);(c_k);\\(\alpha_l):(\beta_m);(\gamma_n);\end{matrix} x,y\right] = \sum_{r,s=0}^{\infty} \frac{\prod_{j=1}^{p}(a_j)_{r+s} \prod_{j=1}^{q}(b_j)_r \prod_{j=1}^{k}(c_j)_s}{\prod_{j=1}^{l}(\alpha_j)_{r+s} \prod_{j=1}^{m}(\beta_j)_r \prod_{j=1}^{n}(\gamma_j)_s} \frac{x^r y^s}{r!s!}, \qquad (1.1)$$

where for convergence

$$\left. \begin{matrix} p+q < l+m+1, p+k < l+n+1, |x|<\infty, |y|<\infty \\ \text{or} \\ p+q < l+m+1, p+k < l+n+1, |x|<\infty, |y|<\infty, \text{and} \\ |x|^{\frac{1}{p-1}} + |y|^{\frac{1}{p-1}} < 1, \text{if } p>l; \max\{|x|,|y|\} < 1, \text{if } p<l \end{matrix} \right\} \qquad (1.2)$$

where $\prod_{j=1}^{p}(a_j)_{r+s} = (a_1)_{r+s}(a_2)_{r+s} \cdot ... \cdot (a_p)_{r+s}$, with similar interpretations for $\prod_{j=1}^{l}(\alpha_j)_{r+s}$, et cetera.

These hypergeometric functions appear in the solution of the partial differential equations which are dealt with harmonic analysis method (see [6]). It is noted that Riemann functions and the fundamental solutions of the degenerate second-order partial differential equations are expressible by means of hypergeometric functions of several variables (see [2,3,4,9-16,20,21]). Therefore, in investigation of boundary-value problems for these partial differential equations, we need to study the solution of the system of hypergeometric functions and find explicit linearly independent solutions (see [10-13]). The function $F_{l:i;j}^{p:q;k}$ contains a large number of Kampé de Fériet type functions.

Here, we choose the functions $F_{2:0;1}^{1:2;1}$ and $F_{1:0;1}^{0:2;1}$ defined, respectively, by the following double series (see [19]):

$$F_{2:0;1}^{1:2;1}\left[\begin{array}{c}a:b,c;d;\\e,f:\ -;g;\end{array}x,y\right]=\sum_{m,n=0}^{\infty}\frac{(a)_{m+n}(b)_m(c)_m(d)_n}{(e)_{m+n}(f)_{m+n}(g)_n}\frac{x^m}{m!}\frac{y^n}{n!}, \quad (1.3)$$

and

$$F_{1:0;1}^{0:2;1}\left[\begin{array}{c}-:b,c;d;\\e:\ -;g;\end{array}x,y\right]=\sum_{m,n=0}^{\infty}\frac{(b)_m(c)_m(d)_n}{(e)_{m+n}(g)_n}\frac{x^m}{m!}\frac{y^n}{n!}, \quad (1.4)$$

where $(a)_n$ denotes the Pochhammer symbol given by $(a)_n = a(a+1)(a+2)...(a+n-1)$, $(a)_0 = 1$, $(a)_n = \Gamma(a+n)/\Gamma(a)$ and $\Gamma$: being the well-known Gamma function, to find the linearly independent solutions of partial differential equations satisfied by these functions. The regions of convergence of the functions $F_{2:0;1}^{1:2;1}$ and $F_{1:0;1}^{0:2;1}$ are given (1.2). It is easy to verify that the functions $F_{2:0;1}^{1:2;1}$ and $F_{1:0;1}^{0:2;1}$ are functions of the fourth and third orders, respectively.

It is not difficult to see that the functions $F_{2:0;1}^{1:2;1}$ and $F_{1:0;1}^{0:2;1}$ are natural generalizations of the well-known hypergeometric Humbert function [6]

$$\Xi_2(b,c;e;x,y)=\sum_{m,n=0}^{\infty}\frac{(b)_m(c)_m}{(e)_{m+n}}\frac{x^m}{m!}\frac{y^n}{n!}.$$

Interest in the functions (1.3) - (1.4) arose also because, for example, the solution of the Cauchy problem for equation

$$\frac{\partial^2 u}{\partial\xi\partial\eta}+\left[\frac{\alpha}{\eta+\xi}+\frac{\beta}{\eta-\xi}\right]\frac{\partial u}{\partial\xi}+\left[\frac{\alpha}{\eta+\xi}-\frac{\beta}{\eta-\xi}\right]\frac{\partial u}{\partial\eta}+\lambda u=0,\ -\frac{1}{2}<\beta\leq\alpha\leq 0$$

with the conditions

$$u(\xi,\xi)=\tau(\xi),\ 0\leq\xi\leq 1,$$

$$\left[2(1-2\beta)\right]^{-2\beta}\lim_{\eta\to\xi}(\eta-\xi)^{2\beta}\left(u_\xi-u_\eta\right)=\nu(\xi),\ 0<\xi<1,$$

is written out via the functions $F_{1:0;1}^{0:2;1}$ and $\Xi_2$ [7]:

$$u(\xi,\eta) = \gamma_1 \frac{(\eta+\xi)^{-\alpha}}{(\eta-\xi)^{1+2\beta}} \int_\xi^\eta \frac{t^\alpha H(\xi,\eta;t;\lambda)}{(\eta-t)^{-\beta}(t-\xi)^{-\beta}} \tau(t)dt$$

$$-\gamma_1 \frac{(\eta+\xi)^{-\alpha}}{(\eta-\xi)^{1+2\beta}} \int_\xi^\eta \frac{(\eta+\xi-2t)t^\alpha}{(\eta-t)^{-\beta}(t-\xi)^{-\beta}} F_{1:0;1}^{0:2;1}\left[\begin{array}{c}-:\alpha,1-\alpha;\ \beta;\\ \beta:\ -;\ 1+\beta;\end{array}\sigma,\rho\right]\tau'(t)dt$$

$$-\gamma_2(\eta+\xi)^{-\alpha}\int_\xi^\eta (\eta-t)^{-\beta}(t-\xi)^{-\beta}t^\alpha \Xi_2(\alpha,1-\alpha;1-\beta;\sigma,\rho)v(t)dt,$$

where

$$H(\xi,\eta;t;\lambda) = 2(1+2\beta)F_{1:0;1}^{0:2;1} - \frac{\alpha}{t}(\eta+\xi-2t)F_{1:0;1}^{0:2;1} - \frac{\partial}{\partial\sigma}F_{1:0;1}^{0:2;1}\cdot\frac{\partial\sigma}{\partial t} + 4\rho\frac{\partial}{\partial\rho}F_{1:0;1}^{0:2;1},$$

$$F_{1:0;1}^{0:2;1} \equiv F_{1:0;1}^{0:2;1}\left[\begin{array}{c}-:\alpha,1-\alpha;\ \beta;\\ \beta:\ -;\ 1+\beta;\end{array}\sigma,\rho\right],\ \sigma = \frac{(\eta-t)(t-\xi)}{2t(\eta+\xi)},\ \rho = \lambda(\eta-t)(t-\xi),$$

$$\gamma_1 = 2^{\alpha-1}\frac{\Gamma(1+2\beta)}{\Gamma^2(1+\beta)},\ \gamma_2 = [2(1-2\beta)]^{2\beta} 2^{\alpha-1}\frac{\Gamma(1-2\beta)}{\Gamma^2(1-\beta)}.$$

According to the theory of multiple hypergeometric functions (see [1]), the system of partial differential equations for the Kampé de Fériet's hypergeometric series of two variables $F_{2:0;1}^{1:2;1}$ and $F_{1:0;1}^{0:2;1}$ are readily seen to be given

$$\left.\begin{array}{c}\left(1+x\frac{\partial}{\partial x}\right)\left(e+x\frac{\partial}{\partial x}+y\frac{\partial}{\partial y}\right)\left(f+x\frac{\partial}{\partial x}+y\frac{\partial}{\partial y}\right)x^{-1}u\\ -\left(a+x\frac{\partial}{\partial x}+y\frac{\partial}{\partial y}\right)\left(b+x\frac{\partial}{\partial x}\right)\left(c+x\frac{\partial}{\partial x}\right)u = 0\\ \left(1+y\frac{\partial}{\partial y}\right)\left(e+x\frac{\partial}{\partial x}+y\frac{\partial}{\partial y}\right)\left(f+x\frac{\partial}{\partial x}+y\frac{\partial}{\partial y}\right)\left(g+y\frac{\partial}{\partial y}\right)y^{-1}u\\ -\left(a+x\frac{\partial}{\partial x}+y\frac{\partial}{\partial y}\right)\left(d+y\frac{\partial}{\partial y}\right)u = 0,\end{array}\right\} \quad (1.5)$$

where $u = F_{2:0;1}^{1:2;1}\left[\begin{array}{c}a:b,c;d;\\ e,f:\ -;\ g;\end{array}x,y\right]$, and

$$\left.\begin{array}{c}\left(1+x\frac{\partial}{\partial x}\right)\left(e+x\frac{\partial}{\partial x}+y\frac{\partial}{\partial y}\right)x^{-1}u - \left(b+x\frac{\partial}{\partial x}\right)\left(c+x\frac{\partial}{\partial x}\right)u = 0\\ \left(1+y\frac{\partial}{\partial y}\right)\left(e+x\frac{\partial}{\partial x}+y\frac{\partial}{\partial y}\right)\left(g+y\frac{\partial}{\partial y}\right)y^{-1}u - \left(d+y\frac{\partial}{\partial y}\right)u = 0,\end{array}\right\} \quad (1.6)$$

where $u = F_{1:0;1}^{0:2;1}\left[\begin{array}{c}-:b,c;d;\\ e:\ -;\ g;\end{array}x,y\right]$, respectively.

## 2. Solving the systems of Partial Differential Equations

Starting from (1.5) and by making use of come elementary calculations, we define the system of third and fourth orders partial differential equations:

$$x^2(1-x)u_{xxx} + xy(2-x)u_{xxy} + y^2 u_{xyy} + x[e+f+1-(a+b+c+3)x]u_{xx} +$$
$$+y[e+f+1-(b+c+1)x]u_{xy} + \{ef - [a(b+c+1)+(b+1)(c+1)]x\}u_x \qquad (2.1)$$
$$-bcyu_y - abcu = 0,$$

$$y^3 u_{yyyy} + x^2 y u_{xxyy} + 2xy^2 u_{xyyy} + (e+f+g+3)y^2 u_{yyy} + (e+f+2g+3)xyu_{xyy}$$
$$+gx^2 u_{xxy} + \{[(e+1)(f+1)+(e+f+1)g]-y\} yu_{yy} \qquad (2.2)$$
$$+[(e+f+1)g-y]xu_{xy} - dxu_x + [efg - (a+d+1)y]u_y - adu = 0.$$

We note that two equations of system (2.1) and (2.2) are simultaneous, because the hypergeometric function $F_{2:0;1}^{1:2;1}[x, y]$ satisfies the system. Now, in order to find the linearly independent solutions of system (2.1) and (2.2), we will search the solutions in the form

$$u = x^\tau y^\nu \omega, \qquad (2.3)$$

where $\omega$ is an unknown function, and $\tau$ and $\nu$ are constants, which are to be determined. Next, substituting $u = x^\tau y^\nu \omega,$ into system (2.1) and (2.2), we have

$$x^2(1-x)\omega_{xxx} + (2-x)xy\omega_{xxy} + y^2 \omega_{xyy} + [(3\tau+2\nu+e+f+1)-(3\tau+\nu+a+b+c+3)x]x\omega_{xx}$$
$$+[(4\tau+2\nu+e+f+1)-(2\tau+b+c+1)x]y\omega_{xy} + \tau x^{-1} y^2 \omega_{yy}$$
$$+\{[\tau(2e+2f+4\nu+3\tau-1)+\nu(e+f+\nu)+ef]$$
$$-[3\tau(\tau-1)+2\tau(\nu+a+b+c+3)+(b+c+1)\nu+(b+1)(a+c+1)+ac]x\}\omega_x \qquad (2.4)$$
$$+\{\tau(2\tau+2\nu+e+f-1)x^{-1} - [\tau(\tau+b+c)+bc]\} y\omega_y$$
$$+\{\tau[(\nu+e-1)(\nu+f-1)+\tau(\tau+2\nu+e+f-2)]x^{-1}$$
$$-[\tau(\tau-1)(\tau+a+b+c+1)+\tau\nu(\tau+b+c)+\tau(b+1)(c+1)+a(b+c+1)\tau+bc(\nu+a)]\}\omega = 0,$$

$$y^3 \omega_{yyyy} + x^2 y \omega_{xxyy} + 2xy^2 \omega_{xyyy} + (e+f+g+2\tau+4\nu+3)y^2 \omega_{yyy}$$
$$+(2\nu+g)x^2 \omega_{xxy} + (e+f+2g+2\tau+6\nu+3)xy\omega_{xyy} + \nu(g+\nu-1)x^2 y^{-1} \omega_{xx}$$
$$+[2\nu(2\tau+3\nu+e+f+2g)+(e+f+2\tau+1)g-y]x\omega_{xy}$$
$$+[(e+1)(f+1)+(e+f+1)g+\tau(\tau+e+f+2g+2)+3\nu(2\tau+2\nu+e+f+g+1)-y]y\omega_{yy} \qquad (2.5)$$
$$+[\nu(\nu+g-1)(2\tau+2\nu+e+f-1)y^{-1} - (\nu+d)]x\omega_x$$
$$+[\nu(\nu-1)(4\nu+6\tau+3e+3f+3g+1)+2\tau\nu(e+f+2g+3)$$
$$+2\nu(e+1)(f+1)+g(e+f+1)(\tau+2\nu)+\tau(\tau-1)(2\nu+g)+efg-(a+d+\tau+2\nu+1)y]\omega_y$$
$$+\nu[(\nu+e-1)(\nu+f-1)(\nu+g-1)+\tau(\nu-1)(\tau+2\nu+e+f+2g-2)+\tau(\tau-1)g+\tau]y^{-1}\omega$$
$$-[d\tau+(\tau+\nu+a+d)\nu+ad]\omega = 0$$

We note that system in (2.4) and (2.5) is analogical to system (2.1) and (2.2), therefore we require that the conditions

$$\begin{cases} \tau = 0, \\ \nu(\nu + g - 1) = 0 \end{cases} \tag{2.6}$$

should be satisfied. It is evident that system (2.6) has two solutions:

$$\tau = 0, \ \nu = 0 \quad \text{and} \quad \tau = 0, \ \nu = 1 - g. \tag{2.7}$$

Finally, substituting both solutions (2.7) into (2.4) and (2.5), we find the two linearly independent solutions of system (2.1) and (2.2):

$$u_1(x, y) = F_{2:0;1}^{1:2;1}\begin{bmatrix} a:b,c;d; \\ e,f: -; g; \end{bmatrix} x, y \end{bmatrix}, \tag{2.8}$$

$$u_2(x, y) = y^{1-g} F_{2:0;1}^{1:2;1}\begin{bmatrix} 1-g+a:b,c;1-g+d; \\ 1-g+e,1-g+f: -; 2-g; \end{bmatrix} x, y \end{bmatrix}. \tag{2.9}$$

The method of establishing the system of partial differential equations (2.4), (2.5) and derivation the solutions (2.8) - (2.9) of this system for the Kampé de Fériet's hypergeometric series of the fourth order and two variables $F_{2:0;1}^{1:2;1}$ detailed above can be applied mutatis mutandis to obtain the systems of partial differential equations and the solutions of the obtained systems for the Kampé de Fériet's hypergeometric series of the third order and two variables $F_{1:0;1}^{0:2;1}$. In this regard, we find each of the following pairs of system of partial differential equations and their solutions, where we put the solutions directly after the corresponding system of partial differential equations:

$$\begin{cases} x(1-x)u_{xx} + yu_{xy} + [e - (b+c+1)x]u_x - bcu = 0, \\ y^2 u_{yyy} + xyu_{xyy} + gxu_{xy} + (e+g+1)yu_{yy} + (eg-y)u_y - du = 0. \end{cases}$$

$$u_1(x, y) = F_{1:0;1}^{0:2;1}\begin{bmatrix} -:b,c;d; \\ e: -; g; \end{bmatrix} x, y \end{bmatrix}$$

and

$$u_2(x, y) = y^{1-g} F_{1:0;1}^{0:2;1}\begin{bmatrix} - :b,c;1+d-g; \\ 1+e-g: -; 2-g; \end{bmatrix} x, y \end{bmatrix}.$$

**References**

[1] P. Appell and J. Kampe de Feriet, *Fonctions Hypergeometriques et Hyperspheriques; Polynomes d'Hermite,* Gauthier - Villars. Paris, 1926.
[2] J. Barros-Neto and I.M. Gelfand, Fundamental solutions for the Tricomi operator, Duke Math.J. 98(3) (1999), 465-483.
[3] J. Barros-Neto and I.M. Gelfand, Fundamental solutions for the Tricomi operator II, Duke Math.J. 111(3) (2002), 561-584.
[4] J. Barros-Neto and I.M. Gelfand, Fundamental solutions for the Tricomi operator III, Duke Math.J. 128(3) (2005), 119-140.

[5] L. Bers, Mathematical Aspects of Subsonic and Transonic Gas Dynamics, Wiley, New York, 1958.
[6] A. Erde'lyi, W. Magnus, F. Oberhettinger and F.G.Tricomi, *Higher Transcendental Functions,* Vol. I, McGraw-Hill Book Company, New York, Toronto and London, 1953.
[7] T.G. Ergashev, The solution of the Cauchy problem for the degenerating hyperbolic equation of second kind, Uzbek Mathematical Journal, 4 (2009), 180-190.
[8] F.I. Frankl, Selected Works in Gas Dynamics, Nauka, Moscow, 1973.
[9] A.J. Fryant, Growth and complete sequences of generalized bi-axially symmetric potentials. *J. Differential Equations*, **31**(1979), 155–164.
[10] A. Hasanov, Fundamental solutions of generalized bi-axially symmetric Helmholtz equation. *Complex Variables and Elliptic Equations*, **52**(2007), 673–683.
[11] A. Hasanov, Some solutions of generalized Rassias's equation, Intern.J.Appl.Math.Stat. 8(M07) (2007), 20-30.
[12] A. Hasanov, The solution of the Cauchy problem for generalized Euler-Poisson-Darboux equation. Intern.J.Appl.Math.Stat. 8(M07) (2007), 30-44.
[13] A. Hasanov, Fundamental solutions for generalized elliptic equation with two perpendicular lines of degeneration. Intern.J.Appl.Math.Stat. 13(8) (2008), 41-49.
[14] A. Hasanov and E.T. Karimov, Fundamental solutions for a class of three-dimensional elliptic equations with singular coefficients. Appl.Math.Letters, 22(2009), 1828-1832.
[16] A. Hasanov, R.B. Seilkhanova and R.D. Seilova, Linearly independent solutions of the system of hypergeometric Exton function $X_{12}$, Contemporary Analysis and Applied Mathematics Vol.3, No.2, 289-292, 2015.
[17] G. Lohofer, Theory of an electromagnetically deviated metal sphere.1:Abcorbed power, SIAM J.Appl.Math.49(1989), 567-581.
[18] A.W. Niukkanen, Generalized hypergeometric series arising in physical and quantum chemical applications, J.Phys.A:Math.Gen. 16 (1983), 1813-1825.
[19] H.M. Srivastava and P.W. Karlsson, *Multipl. Gaussian Hypergeometric Series*, Halsted Press (Ellis Horwood Limited, Chicherster), John Wiley and Sons, New York, Chichester, Brisbane and Toronto, 1985.
[20] R.J. Weinacht, Some properties of generalized axially symmetric Helmholtz potentials. *SIAM J. Math. Anal*. 5(1974), 147-152.
[21] A. Weinstein, Discontinuous integrals and generalized potential theory. *Trans. Amer. Math. Soc*., 63(1948), 342-354.